\documentclass[a4paper,draft]{article}\pdfpagewidth=\paperwidth\pdfpageheight
=\paperheight\usepackage{stix}\usepackage
{titlesec}\usepackage[leqno]{amsmath}\usepackage[hidelinks,final]{hyperref}\usepackage
{amsthm}\usepackage{mathtools}\usepackage{cite}\usepackage{enumitem}\usepackage
{xpatch}\DeclareSymbolFont{stix2int}{LS2}{stix2cal}{m}{n}\DeclareSymbolFont
{stix2op}{LS2}{stix2ex}{m}{n}\makeatletter\stix@MathSymbol\intop\mathop
{stix2int}{"AF}\stix@MathSymbol\iintop\mathop{stix2int}{"B0}\stix@MathSymbol
\iiintop\mathop{stix2int}{"B1}\stix@MathSymbol\sumop\mathop{stix2op}{"B3}\stix@MathSymbol
\prodop\mathop{stix2op}{"B1}\makeatother\titleformat\section{\normalfont
\normalsize\centering\bfseries}{\thesection.}{4pt}{}\titlespacing\section
{0pt}{3ex plus .5ex minus .5ex}{.75ex plus .1ex}\newtheorem{thm}{Theorem}[section]\newtheorem
*{thm*}{Theorem}\newtheorem
*{auxnamedthm*}{\auxthmname}\newtheorem{cor}[thm]{Corollary}\newtheorem
*{cor*}{Corollary}\newtheorem{lem}[thm]{Lemma}\newtheorem*{lem*}{Lemma}\newtheorem
{prop}[thm]{Proposition}\newtheorem*{prop*}{Proposition}\theoremstyle
{definition}\newtheorem*{auxnameddef*}{\auxthmname
}\newtheorem*{rems*}{Remarks}\newtheorem*{exs*}{Exercise}\newtheorem
*{rem*}{Remark}\newtheorem*{defns*}{Definitions}\newtheorem*{defn*}{Definition}
\begin{document}
\par\par\par\par\par\par\par\par\par\par\textbf{\hrule height.4mm\vskip
3mm\centerline{\LARGE Invariant Subspaces for Integration}\centerline
{\LARGE and Differentiation on Convex Domains}\vskip3mm\hrule height.25mm\vskip
2mm\centerline{\large Erik Avelin}\vskip2mm\hrule height.4mm\vskip
3mm} \vspace{.05in} \begin{quotation}\small\noindent{\textbf{Abstract.}}\hskip
1.3mmFor a conformal map $g$ of the unit disc $U$ onto any bounded
convex domain $\Omega$ we determine the invariant subspaces for the
integration operator $T_gf = \int_0^zfg'$ on the Hardy space $H^p$,
by means of a Fock space model. Under some additional regularity assumptions
on $\partial\Omega$ we also treat integration from a boundary point
of $U$, and we determine the residual subspaces (those invariant subspaces
on which the operator has empty spectrum) for the left inverse $D_gf
= f'/g'$ of $T_g$ on the smooth disc algebra $A^\infty$. The latter
result is new even for $\Omega= U$. \end{quotation}\section{Introduction}\label
{1} The aim of this paper is to provide partial answers to the following
two questions. \def\auxthmname{Problem A}\begin{auxnameddef*} Let $U$
be the unit disc, fix $a \in\bar{U}$, let $g$ be holomorphic in $U$,
and suppose the Pommerenke operator \begin{equation}\label{2} (T_gf)(z)
= \int_a^z f(\zeta)g'(\zeta)\;\textit d\zeta\end{equation} is well
defined and bounded on the Hardy space $H^p$. What are its invariant
subspaces? \end{auxnameddef*}\def\auxthmname{Problem B}\begin{auxnameddef*}
Let $A^\infty= \{f\colon D^nf \in A\text{ for }n = 0,1,2,\ldots\}$,
where $(Df)(z) = f'(z)$ and $A$ is the disc algebra. The topology of
$A^\infty$ is that of uniform convergence in all derivatives. What
are the invariant subspaces for $D$ on $A^\infty$? \end{auxnameddef*}By
an \emph{invariant subspace} for a continuous linear operator $S$ on
a topological vector space $X$, we mean a closed subspace $Y \subset
X$ such that $SY \subset Y$ and $Y \ne0$, $Y \ne X$. \par When $a \in
T$ (the unit circle), the integral \eqref{2} is to be interpreted as
the value at $z$ of the unique antiderivative of $fg'$ which has limit
$0$ at $a$. This makes sense for example if $fg' \in H^1$, for then
$\int_0^zfg' \in A$ by Hardy's inequality. \par The operators $T_g$
of Problem~A have attracted much interest in the recent decades, motivated
originally by Pommerenke's observation \cite{19} that they lead to
a very simple proof of the analytic John-Nirenberg inequality. Necessary
and sufficient conditions for their boundedness and compactness have
been determined on a wide range of spaces, including those of Hardy
\cite{4,9,19,21}, Bergman \cite{10}, and Dirichlet \cite{14}. The spectrum
of $T_g$ is also quite well understood, at least on Hardy and Bergman
spaces \cite{5,8,15}. In contrast, almost nothing is known about invariant
subspaces beyond the simplest case of $g(z) = z$, which was solved
by Aleman and Korenblum \cite{7} in 2008. \par Problem~B was posed
to me by Alex Bergman. Its real-variable analogue, concerning ordinary
differentiation on $C^\infty(0,1)$, has been studied by several authors.
Aleman and Korenblum \cite{6} have determined the \emph{residual subspaces},
that is, those invariant subspaces on which the operator has empty
spectrum; the corresponding spectral synthesis problem was subsequently
studied by Aleman, Baranov, and Belov \cite{1} and later solved very
completely by Baranov and Belov \cite{11}. \par Recently, Aleman and
Bergman \cite{3} have for a large class of pairs of operators $(D,
V)$ where $V$ is compact with $\sigma(V) = \{0\}$ and $DV = I$ obtained
a description of the invariant subspaces for $V$ and the residual subspaces
for $D$, generalizing both the result of \cite{6} and the classical
characterization of the Volterra-invariant subspaces of $L^{\mkern
-2.6mu2}$. The operators of Problems~A and~B are interesting examples
of ones to which the techniques of \cite{3} are difficult to apply,
as they rely heavily on the existence of an orthonormal basis among
the eigenvectors $e_\lambda$ of $D$. \par The general trick which is
key to the arguments of both \cite{3}, \cite{6}, and \cite{7}, and
which will also be the foundation of our approaches to Problems~A and~B,
is the conversion of questions of $V\!$-invariance and $D$-residuality
into questions of near invariance (the permission to divide out zeros)
in a space $X$ of entire functions; $X$ is the range of the generalized
Fourier transform $\hat x(\lambda) = \langle e_\lambda,x\rangle$ acting
on the dual of the original space. In the case of Problem~A this transform
is injective only if $\{e^{\lambda g}\colon\lambda\in C\mkern-12mu\mathrlap
/\mkern12mu\}$ spans a dense subspace of $H^p$, and we must therefore
restrict our attention to univalent $g$. \par The hypotheses of \cite
{3} are chosen in such a way that $X$ becomes a de~Branges space. Near
invariance is well understood in these spaces, thanks to de~Branges'
celebrated ordering theorem. Our Problem~B instead gives rise to a
Fock-type space, and the same goes for Problem~A if $g(U)$ is assumed
to be bounded and convex. While near invariance is in general quite
poorly understood in Fock-type spaces, a recent result on backward-shift
invariance due to Aleman, Baranov, Belov, and Hedenmalm \cite{2} has
allowed us to prove the following. \def\auxthmname{Theorem A.1}\begin{auxnamedthm*}
Suppose $0 < p < \infty$ and $g$ is a conformal map of $U$ onto a bounded
convex domain. Fix $a \in U$ and put \begin{equation*} B_a(z) = \frac
{z-a}{1-\bar az}. \end{equation*} Then \begin{enumerate}[label={\textup
{(\textit{\alph*})}},listparindent=\parindent,parsep=0pt]\item$T_g$
is a compact operator $H^p \to H^p$ with $\sigma(T_g) = \{0\}$. \item
$B_a^nH^p$ is $T_g$-invariant for $n=1,2,3,\ldots$. \item\label{3}
No other subspace of $H^p$ is $T_g$-invariant. \end{enumerate}\end{auxnamedthm*}It
is obvious that \ref{3} is false for any $g$ which is not conformal;
consider the closure of the polynomials in $g$. \par For Problem~B
and for the case $a \in T$ of Problem~A we will instead resort to a
weak but very general result about near invariance: Hadamard's factorization
theorem. The technique, which was also used in \cite{7}, is the following.
\par Suppose $X$ is any topological vector space of entire functions
of exponential type. A closed subspace $N \subset X$ ($N \ne0$, $N
\ne X$) is said to be \emph{nearly invariant} if $F(z)/(z-\lambda)
\in N$ whenever $F \in N$ and $F(\lambda) = 0$. (Thus we do not permit
$N$ to have common zeros, unlike in some literature.) For any $F \in
X$, logarithmic differentiation of the Hadamard expansion \begin{equation}\label
{4} F(z) = e^{b+cz}z^m\prod_{n=1}^\infty(1-z/\lambda_n)e^{z/\lambda
_n} \end{equation} leads to a locally uniformly convergent series \begin{equation}\label
{5} F'(z) = cF(z) + m\frac{F(z)}{z} + \sum_{n=1}^\infty\frac{zF(z)}{\lambda
_n(z-\lambda_n)}. \end{equation} If $F \in N$ and $N$ is nearly invariant,
all terms of \eqref{5} belong to $N$. Thus, we can draw the useful
conclusion that $F' \in N$, provided that the series converges in the
topology of $X$. In the case of Problem~A this convergence turns out
to be related to the behavior of the quantity \begin{equation}\label
{6} \textit d_g(\lambda) = \inf\{\Vert R_\lambda-S\Vert\colon\mathscr
R(S) \subset C\mkern-12mu\mathrlap/\mkern12mue^{\lambda g}\} \end{equation}
for large $|\lambda|$. The norm in \eqref{6} is the $H^p$ operator
norm, $\mathscr R(S)$ denotes the range of the operator $S$, and \begin{equation}
R_\lambda= T_g(I-\lambda T_g)^{-1}. \end{equation} \par Since this
technique does not rely on $X$ being a Fock-type space it could potentially
be used to extend Theorem~A.1 to more general conformal maps $g$, but
the author has been unable to do this --- the method of estimating
$\textit d_g(\lambda)$ presented in Section~\ref{73} breaks down completely
if $g(U)$ is not convex. We do however obtain the following two results.
\def\auxthmname{Theorem A.2}\begin{auxnamedthm*} Suppose $1 \le p <
\infty$ and $g$ is a conformal map of $U$ onto a Jordan domain whose
boundary is $C^2$ with positive curvature. Fix $a \in T$ and put \begin{equation*}
S_a^t(z) = \exp\bigg\{-t\frac{a+z}{a-z}\bigg\}. \end{equation*} Then
\begin{enumerate}[label={\textup{(\textit{\alph*})}},listparindent=\parindent
,parsep=0pt]\item$T_g$ is a compact operator $H^p \to H^p$ with $\sigma
(T_g) = \{0\}$. \item$S_a^tH^p$ is $T_g$-invariant for any $t > 0$.
\item\label{7} No other subspace of $H^p$ is $T_g$-invariant. \end{enumerate}\end{auxnamedthm*}Explicitly,
the boundary curve is to have a parametrization $\Gamma\in C^2$ with
$|\Gamma'| = 1$ and $|\Gamma''| > 0$. Note that this implies convexity.
\par By the spectrum $\sigma(\Theta)$ of an inner function $\Theta
= cBS$, $|c| = 1$, we will mean the closure of $Z(B) \cup\operatorname
{supp}\mu$, where $Z(B)$ is the zero set of the Blaschke product $B$
and \begin{equation*} S(z) = \exp\bigg\{-\int_{T}\frac{w+z}{w-z}\;\textit
d\mu(w)\bigg\}. \end{equation*} \par A compact set $K \subset T$ of
Lebesgue measure $0$ is called a \emph{Carleson set} if the lengths
$\delta_i$ of the arcs forming the complement of $K$ satisfy $\sum
\delta_i\log\delta_i > -\infty$. It is well known (see e.g.\ \cite
{22}) that $K$ is a Carleson set if and only if the space \begin{equation*}
A_K^\infty= \{f \in A^\infty\colon D^nf = 0\text{ on }K\text{ for }n
= 0,1,2,\ldots\} \end{equation*} contains some $f \ne0$. \par\def\auxthmname
{Theorem B}\begin{auxnamedthm*} Let $g$ be a conformal map of $U$ onto
a Jordan domain whose boundary is $C^\infty$ with positive curvature.
Write $(D_gf)(z) = f'(z)/g'(z)$. Then \begin{enumerate}[label={\textup
{(\textit{\alph*})}},listparindent=\parindent,parsep=0pt]\item$D_g$
is a bounded operator $A^\infty\to A^\infty$. \item\label{8} $SA_K^\infty
$ is $D_g$-invariant for any Carleson set $K$ and any singular inner
function $S$ with $\sigma(S) \subset K$. \item\label{9} $S_a^tA_a^\infty
$ is $D_g$-residual for any $a \in T$ and any $t \ge0$. \item\label
{10} No other subspace of $A^\infty$ is $D_g$-residual. \end{enumerate}\end{auxnamedthm*}Of
course, the eigenvectors $e^{\lambda g}$ of $D_g$ give rise to invariant
subspaces which are not of the form $SA_K^\infty$. We conjecture that
every $D_g$-invariant subspace containing no $e^{\lambda g}$ is as
in \ref{8}. By \cite{22} this is equivalent to the statement that every
such invariant subspace is an ideal of the algebra $A^\infty$. \par
Let $A(\Omega)$ consist of all functions which are holomorphic in $\Omega
$ and continuous on $\bar\Omega$. Define $A^\infty(\Omega)$ as in Problem~B,
with $A(\Omega)$ in place of $A$. If $g$ is as in Theorem~B and $\Omega
= g(U)$, then the composition $C_gf = f\circ g$ is a bounded linear
isomorphism of $A^\infty(\Omega)$ onto $A^\infty$ with $D_g = C_gDC_g^{-1}$.
Thus, Theorem~B completely describes the residual subspaces for differentiation
on $A^\infty(\Omega)$. Similarly, Theorems~A.1 and~A.2 describe the
invariant subspaces for integration on the conformally invariant Hardy
space $H^p(\Omega)$. \par The layout of the paper is as follows. Section~\ref
{11} contains the verifications of parts~(\textit{a}) and~(\textit
{b}) of all three theorems and part~\ref{9} of Theorem~B. Part~\ref
{3} of Theorem~A.1 is proved in Section~\ref{28}. In Section~\ref{48}
we present in parallel the arguments for \ref{7} of Theorem~A.2 and
\ref{10} of Theorem~B assuming only that $g(U)$ is Jordan with smooth
enough boundary and that $\textit d_g(\lambda) = O(|\lambda|^{-\delta
})$ for some $\delta> 0$. We also show that Theorem~A.1 is true under
these hypotheses, regardless of whether $g(U)$ is convex. The proof
that our assumption of positive curvature gives control of $\textit
d_g(\lambda)$ is in Section~\ref{73}. \par I am indebted to Alexandru
Aleman and Alex Bergman for many useful suggestions, and to Adem Limani
for a number of helpful comments on the manuscript. \section{Examples
of Invariant Subspaces}\label{11} It is clear that any conformal map
of $U$ onto a domain of finite area belongs to the Dirichlet space
and therefore to $\operatorname{VMOA}$. By the Kellogg-Warschawski
theorem (see Theorems~3.5 and~3.6 of \cite{20}), the hypotheses of
Theorems~A.2 and~B imply respectively that $\log g' \in A$ and $\log
g' \in A^\infty$. The following propositions therefore establish the
easy halves of Theorems~A.1, A.2, and~B. \begin{prop}\label{12} Let
$0 < p < \infty$ and fix $a \in U$. Then \begin{enumerate}[label={\textup
{(\textit{\alph*})}},listparindent=\parindent,parsep=0pt]\item\label
{13} $T_g$ is bounded on $H^p$ if and only if $g \in\operatorname{BMOA}$.
\item\label{14} $T_g$ is compact on $H^p$ if and only if $g \in\operatorname
{VMOA}$. \item\label{15} When $g \in\operatorname{BMOA}$, $B_a^nH^p$
is $T_g$-invariant for $n = 1,2,3,\ldots$. \item\label{16} When $g
\in\operatorname{VMOA}$, $\sigma(T_g) = \{0\}$. \end{enumerate}\end{prop}
\begin{proof}[\textsc{proof}] For \ref{13} and \ref{14}, see \cite
{4}. A function $f \in H^p$ is divisible by $B_a^n$ if and only if
$f$ has at least $n$ zeros at $a$. Clearly, $T_gf$ has at least one
more zero at $a$ than $f$ has, and this gives \ref{15}. Part~\ref{16}
is clear from \ref{14}, since $T_g$ has no point spectrum. (More generally,
Limani and Malman \cite{15} have proved that $\sigma(T_g) = \{0\}$
if and only if $g$ is approximable by $H^\infty$-functions in the $\operatorname
{BMO}$-norm.) \end{proof} \begin{prop}\label{17} Let $1 \le p < \infty
$, fix $a \in\bar{U}$ and let $g' \in H^q$, $1/p + 1/q = 1$. Then \begin{enumerate}[label={\textup
{(\textit{\alph*})}},listparindent=\parindent,parsep=0pt]\item\label
{18} $T_gH^p \subset A$. \item\label{19} $T_g$ is compact on $H^p$,
and $\sigma(T_g) = \{0\}$. \item\label{20} When $a \in T$, $S_a^tH^p$
is $T_g$-invariant for $t > 0$. \end{enumerate}\end{prop} \begin{proof}[\textsc
{proof}] If $f(z) = \sum c_nz^n$, Hardy's inequality (applied to $f'$)
says that \begin{equation}\label{21} \sum|c_n| \le|f(0)| + \pi\Vert
f'\Vert_1, \end{equation} where the norm is that of $H^1$. This immediately
yields \ref{18} and shows that $g \in A \subset\operatorname{VMOA}$.
Because $\int_a^z = \int_0^z - \int_0^a$ on $H^1$, \ref{19} follows
from Proposition~\ref{12}. For \ref{20} we appeal to Lemma~3.1 of \cite
{7} which says that a function $f \in H^p$ is divisible by $S_a^t$
if and only if \begin{equation}\label{22} \limsup_{r \to1-}\,(1-r)\log
|f(ra)| \le-2t. \end{equation} If $f \in S_a^tH^p$, then $fg' \in S_a^tH^1$.
Therefore, $fg'$ satisfies \eqref{22}, and it is straightforward to
verify that so does $\int_a^zfg'$. \end{proof} \begin{lem}\label{23}
If $S$ is a singular inner function, then $SA_{\sigma(S)}^\infty= A^\infty
\cap SH^1$ is $D$-invariant in $A^\infty$. \end{lem} \begin{proof}[\textsc
{proof}] That $SA_{\sigma(S)}^\infty$ is closed and equals $A^\infty
\cap SH^1$ is part of Theorem~4.7 of \cite{22}. The differentiation
invariance follows from the theorem of Caughran \cite{12} which says
that $f' \in SH^1$ whenever $f' \in H^1$ and $f \in SH^1$. \end{proof}
\begin{prop}\label{24} If $\log g' \in A^\infty$, then \begin{enumerate}[label={\textup
{(\textit{\alph*})}},listparindent=\parindent,parsep=0pt]\item\label
{25} $D_g$ is bounded on $A^\infty$. \item\label{26} $SA_K^\infty$
is $D_g$-invariant for any Carleson set $K$ and any singular inner
function $S$ with $\sigma(S) \subset K$. \item\label{27} $S_a^tA_a^\infty
$ is $D_g$-residual for $a \in T$ and $t \ge0$. \end{enumerate}\end{prop}
\begin{proof}[\textsc{proof}] \ref{25} This is trivial. \par\ref{26}
Since every $f \in A \cap SH^1$ vanishes on $\sigma(S)$, Lemma~\ref
{23} implies that \begin{equation*} SA_K^\infty\subset SA_{\sigma(S)}^\infty
\subset A_{\sigma(S)}^\infty. \end{equation*} Since $S$ is holomorphic
off $\sigma(S)$, it follows easily that $SA_K^\infty= SA_{\sigma(S)}^\infty
\cap A_K^\infty$, an intersection of $D_g$-invariant subspaces. \par
\ref{27} By Proposition~\ref{17}, the operator $R_\lambda= T_g(I-\lambda
T_g)^{-1}$ maps $H^1$ into $A$ and $S_a^tH^1$ into $S_a^tH^1$. Since
$D_gR_\lambda= I + \lambda R_\lambda$ it follows via Lemma~\ref{23}
that $R_\lambda$ maps $A^\infty$ into $A^\infty$ and $S_a^tA_a^\infty
$ into $S_a^tA_a^\infty$. By the closed graph theorem, $R_\lambda$
is bounded on $A^\infty$. \par As we just noted, $(D_g-\lambda I)R_\lambda
= I$. On $S_a^tA_a^\infty$ we also have $T_gD_g = I$ and therefore
$R_\lambda(D_g-\lambda I) = I$. Thus, $D_g$ restricted to $S_a^tA_a^\infty
$ has empty spectrum. \end{proof} \section{Integration from an Interior
Point} \label{28} Fix $a$ and $g$ as in Theorem~A.1. Let $g(a) = 0$,
$\Omega= g(U)$, and $\varphi= 1/\sqrt{g'}$. It involves no loss of
generality to assume that $a = 0$, since the composition $C_af = f
\circ B_a$ is a bounded linear automorphism of $H^p$ with \begin{equation*}
(C_a^{-1}T_gC_af)(z) = \int_0^z fg_a',\qquad g_a = g \circ B_a^{-1}.
\end{equation*} \par We will begin by showing that Theorem~A.1 is true
if $H^p$ is replaced by the space $\varphi H^2$, normed so that $\Vert
\varphi f\Vert$ equals the $H^2$-norm of $f$. This space is the isometric
image of the Smirnov space $E^2(\Omega)$ under $f \to f \circ g$. See
Chapter~10 of \cite{13}. \begin{prop}\label{29} If $g$ and $\varphi
H^2$ are as above, then \begin{enumerate}[label={\textup{(\textit{\alph
*})}},listparindent=\parindent,parsep=0pt]\item\label{30} $1/g' \in
H^\infty$ and $g' \in H^p$ for some $p > 1$. \item\label{31} $T_g$
is compact on $\varphi H^2$ with $\sigma(T_g) = \{0\}$. \item\label
{32} There exist $q < \infty$ such that $H^q \subset\varphi H^2 \subset
H^2$. \item\label{33} If $p > 0$ and $q < \infty$, then $T_g^nH^p \subset
H^q$ for some positive integer $n$. \item\label{34} The polynomials
are dense in $\varphi H^2$. \end{enumerate}\end{prop} \begin{proof}[\textsc
{proof}] Recall that every bounded convex domain is Jordan and has
rectifiable boundary. In particular, $g' \in H^1$ and $g \in A$. By
Corollary~3.16 of \cite{20}, \begin{equation}\label{35} g'(z) = g'(0)\exp
\bigg\{{-}\frac{1}{\pi}\int_0^{2\pi}\log(1-e^{-it}z)\;\textit d\beta
(t)\bigg\} \end{equation} where $\beta$ is an increasing right-continuous
function with $\beta(t+2\pi) = \beta(t) + 2\pi$. Geometrically, $e^{i\beta
(t)}$ is the direction of the forward half-tangent of $\partial\Omega
$ at $g(e^{it})$. The first half of \ref{30} is clear from \eqref{35}.
\par A Jordan domain $\Omega$ with rectifiable boundary is called a
\emph{Lavrentiev domain} if there exists a constant $M < \infty$ such
that \begin{equation*} \Lambda(J(z,w)) \le M|z-w|\,\qquad(z,w \in\partial
\Omega) \end{equation*} where $J(z,w)$ is the shorter arc of $\partial
\Omega$ connecting $z$ to $w$ and $\Lambda(J)$ denotes the length of
$J$. A quick way to see that bounded convex domains are Lavrentiev
(pointed out to me by Eskil Rydhe) is the following. There exists $\epsilon
> 0$ such that $\beta(s) - \beta(t) \le\pi- \epsilon$ whenever $t <
s < t + \epsilon$. If $z = g(e^{it})$ and $w = g(e^{is})$ it follows
that the forward half-tangent at $z$ meets at some point $\alpha$ the
backward half-tangent at $w$. The angle at $\alpha$ of the triangle
$[z,\alpha,w]$ is at least $\epsilon$, and the arc $J$ from $z$ to
$w$ lies in $[z,\alpha,w]$. Thus, for some $M$ depending only on $\epsilon
$, \begin{equation*} \Lambda(J) \le2\pi\operatorname{diam}J \le M|z-w|.
\end{equation*} \par For any conformal map $g$ of $U$ onto a Lavrentiev
domain, Theorem~7.11 of \cite{20} shows that $g'$ is outer with $|g'|$
satisfying the Coifman-Fefferman $A_\infty$ condition. We shall require
only two consequences of the $A_\infty$ condition, both of which are
proved in \cite{8}. The first is the classical fact that the second
half of \ref{30} holds. The second is that the operator \begin{equation*}
(Sf)(z) = \frac{1}{\varphi(z)}\int_0^z\varphi f' \end{equation*} is
bounded on $H^2$. Since $T_g$ is compact on $H^2$, it follows that
$ST_g$ is compact on $H^2$. This implies \ref{31}, because $T_g(\varphi
f) = \varphi ST_gf$. \par Parts~\ref{32}, \ref{33}, and~\ref{34} are
corollaries of \ref{30}, by Hölder's inequality, Theorem~1 of \cite
{4}, and Beurling's theorem, respectively. \end{proof} If $\mu$ is
a positive measure on $C\mkern-12mu\mathrlap/\mkern12mu$ we will write
$A^2(\mu)$ for the space of all entire functions in $L^{\mkern-2.6mu2}(\mu
)$. By a theorem of Lutsenko and Yulmukhametov \cite{17} (see also
\cite{16} and \cite{18}), our assumption that $\Omega$ is bounded and
convex implies that the transform \begin{equation}\label{36} (\Psi
h)(w) = \langle e^{wg}, h\rangle\end{equation} takes the dual space
$E$ of $\varphi H^2$ onto some $A^2(\nu)$, with \begin{equation}\label
{37} \frac{1}{C}\Vert h\Vert^2 \le\int|\Psi h|^2\;\textit d\nu\le C\Vert
h\Vert^2\,\qquad(h \in E) \end{equation} for some constant $C < \infty
$. The fact that there exists a conjugate-linear isometry between $E$
and $\varphi H^2$ will be of no relevance in the following discussion,
but we will use the fact that norm-closed subspaces of $E$ are weak*-closed.
\par For $\lambda\in C\mkern-12mu\mathrlap/\mkern12mu$ and entire functions
$F$, define \begin{equation*} (L_\lambda F)(w) = \frac{F(w)-F(\lambda
)}{w-\lambda}. \end{equation*} Let $L = L_0$. Because of the normalization
$g(a) = 0$, direct computation shows that $\Psi T_g^* = L\Psi$ on $E$.
By \ref{31} of Proposition~\ref{29} it follows that $L$ is a compact
operator on $A^2(\nu)$ with $\sigma(L) = \{0\}$. One checks easily
that $L_\lambda= L(I-\lambda L)^{-1}$. \par A subspace of $A^2(\mu
)$ is called \emph{backward-shift invariant} if it is an invariant
subspace for every $L_\lambda$. The preceding paragraph shows that
$Y$ is $T_g$-invariant in $\varphi H^2$ if and only if $\Psi Y^\perp
$ is backward-shift invariant in $A^2(\nu)$. We will appeal to the
following theorem from Section~2 of \cite{2}. \begin{thm}\label{38}
If $\mu$ is a positive measure such that \begin{enumerate}[label={\textup
{(\textit{\roman*})}},listparindent=\parindent,parsep=0pt]\item\label
{39} convergence in $A^2(\mu)$ implies locally uniform convergence,
\item\label{40} every $f \in A^2(\mu)$ is of finite order, and \item
\label{41} $A^2(\mu)$ contains the polynomials as a dense subspace,
\end{enumerate}then every backward-shift invariant subspace of $A^2(\mu
)$ is some $\mathscr P_n$, the space of all polynomials of degree at
most $n$. \end{thm} It is clear from \eqref{36} and \eqref{37} that
\ref{39} and \ref{40} hold for $\mu= \nu$. If $m$ denotes normalized
Lebesgue measure on the unit circle and $h_n$ is the functional defined
on $\varphi H^2$ by \begin{equation*} \langle f, h_n\rangle= \int f(z)\bar
z^n\;\textit dm(z), \end{equation*} then $h_n$ is bounded by \ref{32}
of Proposition~\ref{29}, and we have \begin{equation*} (\Psi h_n)(w)
= \sum_{k=0}^\infty\frac{w^k}{k!}\int g^k(z)\bar z^n\;\textit dm(z).
\end{equation*} Since $g^k$ has $k$ zeros at $0$, it follows that $\Psi
$ maps $\mathscr Q_n = \operatorname{span}\{h_0,\ldots,h_n\}$ into
$\mathscr P_n$. Since $\Psi$ is one-to-one and $\dim\mathscr Q_n =
\dim\mathscr P_n < \infty$ we must have $\Psi\mathscr Q_n = \mathscr
P_n$. By the Hahn-Banach theorem the span of all $h_n$ is dense in
$E$, and \ref{41} follows.\par By Theorem~\ref{38}, $T_g$ is \emph
{unicellular} on $\varphi H^2$, that is, the invariant subspaces are
totally ordered by set inclusion. \par If $T$ is a continuous linear
operator on a topological vector space $X$, we will say that a chain
of invariant subspaces $X_1 \supset X_2 \supset\cdots$ is \emph{complete}
if no nonzero $x \in X$ lies in every $X_n$ and if $\dim X_{n-1}/X_n
= 1$ for $n = 1,2,3,\ldots$, where $X_0 = X$. \par Since Proposition~\ref
{12} yields a complete chain of $T_g$-invariant subspaces of $H^p$,
Theorem~A.1 will be proved as soon as we show that $T_g$ is unicellular
on $H^p$. We will use an elementary operator-theoretic trick to transfer
the unicellularity from $\varphi H^2$. \par\begin{prop}\label{42} Let
$T$ be a continuous linear operator on a topological vector space $Z$.
Let $X$ and $Y$ be F-spaces continuously contained in $Z$, with $TX
\subset X$, $TY \subset Y$. Let $X_n$ be the $X$-closure of $T^nX$
and define $Y_n$ similarly for $Y$. Suppose \begin{enumerate}[label={\textup
{(\textit{\roman*})}},listparindent=\parindent,parsep=0pt]\item\label
{43} $X_1 \supset X_2 \supset\cdots$ is a complete chain of $T$-invariant
subspaces of $X$, \item\label{44} $Y_1 \supset Y_2 \supset\cdots$ is
a complete chain of $T$-invariant subspaces of $Y$, \item\label{45}
$T^kX \subset Y$ and $T^kY \subset X$ for some positive integer $k$.
\end{enumerate}Then $T$ is unicellular on $X$ if and only if $T$ is
unicellular on $Y$. \end{prop} \begin{proof}[\textsc{proof}] The closed
graph theorem guarantees that $T^k$ maps $X$ continuously into $Y$
and vice versa, and that $T$ is a bounded operator on $X$ and on $Y$.
\par Fix $x_0$ so that $X = C\mkern-12mu\mathrlap/\mkern12mux_0 \oplus
X_1$. Since $TX_1 \subset X_2$, we must have $Tx_0 \notin X_2$, for
otherwise $TX \subset X_2$. Writing $x_n = T^nx_0$, it follows that
\begin{equation*} X = C\mkern-12mu\mathrlap/\mkern12mux_0 \oplus\cdots
\oplus C\mkern-12mu\mathrlap/\mkern12mux_{n-1} \oplus X_n\,\qquad(n
= 1,2,3,\ldots). \end{equation*} In particular, $T$ is injective on
$X$. \par Since $X_n$ is $T$-invariant, $T$ is well defined on $X/X_n$.
The vectors $x_0,\ldots,x_{n-1}$ project to a basis of $X/X_n$ in which
the the matrix of $T$ is a single (transposed) Jordan block. Thus $T$
is unicellular on $X/X_n$. To show that $T$ is unicellular on $X$ it
is therefore enough to show that every $T$-invariant subspace of $X$
contains an $X_n$. \par Assume $T$ is unicellular on $Y$ and let $M$
be $T$-invariant in $X$. By \ref{44}, the $Y$-closure of $T^kM$ equals
$Y_n$ for some $n$. Since $T^k$ maps $Y$ continuously into $X$, it
follows that $T^{2k}M$ is $X$-dense in $T^kY_n$. In particular, $T^kY_n
\subset M$. This implies that $X_{n+2k} \subset M$, because $T^{n+2k}X
\subset T^{n+k}Y \subset T^kY_n$. \end{proof} Take $Z = H(U)$, $X =
H^p$, $Y = \varphi H^2$, and $T = T_g$, where $H(U)$ is the space of
all holomorphic functions in $U$ with the topology of locally uniform
convergence. It is clear from Propositions~\ref{12} and~\ref{29} that
all hypotheses of Proposition~\ref{42} are satisfied, except perhaps
\ref{43} and \ref{44}. \par If $n$ is a positive integer and $P$ is
a polynomial, \ref{30} of Proposition~\ref{29} shows that $B_a^nP =
T_gf$ for some $f \in B_a^{n-1}H^\infty$. Therefore, the $H^p$-closure
of $T_g(B_a^{n-1}H^p)$ is $B_a^nH^p$. It follows by induction that
$X_n = B_a^nH^p$, and this gives \ref{43}. The same argument applies
to $\varphi H^2$ in place of $H^p$, by \ref{34} of Proposition~\ref
{29}. This completes the proof of Theorem~A.1. \par For $0 < \alpha
\le1$, let $\Lambda_\alpha$ be the space of all $f \in H(U)$ for which
$(1-|z|)^{1-\alpha}|f'(z)|$ is bounded. By Theorem~1 of \cite{4}, \ref
{33} of Proposition~\ref{29} is true for any holomorphic function $g$
which belongs to some $\Lambda_\alpha$. As a by-product of the preceding
argument we therefore obtain the following result. \begin{prop}\label
{46} Suppose $a \in U$, $g \in\Lambda_\alpha$, and $1/g' \in H^p$ for
some $\alpha> 0$ and all $p < \infty$. For $0 < p < \infty$, let $Q(p)$
be the statement that every $T_g$-invariant subspace of $H^p$ is of
the form $B_a^nH^p$. Then $Q(p)$ is true for some $p$ only if it is
true for all $p$. \end{prop} \par For the sake of completeness we also
give an analogue for $a \in T$. \par\begin{prop}\label{47} Suppose
$a \in T$, $g' \in H^\infty$, and $1/g' \in H^p$ for all $p < \infty
$. For $1 \le p < \infty$, let $Q(p)$ be the statement that every $T_g$-invariant
subspace of $H^p$ is of the form $S_a^tH^p$. Then $Q(p)$ is true for
some $p$ only if it is true for all $p$. \end{prop} \begin{proof}[\textsc
{proof}] Suppose $Q(p)$ holds, let $1 \le q < \infty$, and let $M \subset
H^q$ be $T_g$-invariant. By Proposition~\ref{17}, $T_gH^q \subset H^p$.
By hypothesis, the $H^p$-closure of $T_gM$ is some $S_a^tH^p$. In particular,
(\textit{i}) $T_gM \subset S_a^tH^p$ and (\textit{ii}) $T_g(S_a^tH^p)
\subset M$. By the theorem of Caughran stated in the proof of Lemma~\ref
{23}, (\textit{i}) implies that $M \subset S_a^tH^q$. The derivative
of $S_a^t$ is $\psi S_a^t$, where $(z-a)^2\psi(z)$ is constant. Therefore,
(\textit{ii}) shows that $M$ contains $S_a^tP$ for any polynomial $P$
with a double zero at $a$. Since $a \in T$, this makes $S_a^tH^q \subset
M$. \end{proof} \par Of course, the above techniques can also be applied
to spaces other than $H^p$, as long as $g$ is sufficiently well behaved.
For instance, if $a \in U$ and $g$ is a conformal map of $U$ onto a
convex domain with $g' \in H^\infty$, it is easy to show using Proposition~\ref
{42} that the $T_g$-invariant subspaces of the Bergman space $A^2 =
H(U) \cap L^{\mkern-2.6mu2}(A)$, where $A$ denotes area measure, are
precisely $B_a^nA^2$ for $n = 1,2,3,\ldots$. See also Section~4 of
\cite{7}. \section{Integration from a Boundary Point and Differentiation}\label
{48} In this section, $1 < p < \infty$. We fix $a \in\bar{U}$ and a
conformal map $g$ of $U$ onto a Jordan domain $\Omega$. The standing
assumptions on $g$ are \begin{enumerate}[label={\textup{(\textit{\roman
*})}},listparindent=\parindent,parsep=0pt]\item\label{49} $\textit
d_g(\lambda) = O(|\lambda|^{-\delta})$ as $|\lambda| \to\infty$ for
some $\delta> 0$, and \item\label{50} $\log g' \in A$ or $\log g' \in
A^\infty$, depending on the context. \end{enumerate}The quantity $\textit
d_g(\lambda)$ was defined in \eqref{6} and is the infimum of the norms
of all operators on $H^p$ which satisfy one (hence all) of the three
conditions listed in Proposition~\ref{52} below. In particular, $\textit
d_g(\lambda)$ does not depend on $a$. Results about operators on $H^p$
will be presented in parallel with results about operators on $A^\infty
$. The stronger version of \ref{50} is needed only for the $A^\infty
$ results. We will prove in Section~\ref{73} that \ref{49} is true
under the hypotheses of Theorem~A.2. \par The dual spaces of $H^p$
and $A^\infty$ will be denoted by $H_*^p$ and $A_*^\infty$, respectively.
The usual identifications of $H_*^p$ and $A_*^\infty$ with spaces of
functions or distributions will not be needed, but we will use the
fact that norm-closed subspaces of $H_*^p$ are weak*-closed. \par For
$f \in A^\infty$ and $n = 0,1,2,\ldots$, Hardy's inequality~\eqref
{21} implies that \begin{equation*} \Vert D^nf\Vert_p \le\Vert D^nf\Vert
_\infty\le\Vert D^nf\Vert_p + \pi\Vert D^{n+1}f\Vert_p. \end{equation*}
Every $D^n$ is a polynomial in $D_g$ with coefficients in $A^\infty
$. The same is true if $D$ and $D_g$ are interchanged. Therefore, the
topology defined on $A^\infty$ in Section~\ref{1} is compatible with
the seminorms $\Vert D_g^nf\Vert_p$. If $u \in A_*^\infty$ there thus
exist $h_0,\ldots,h_n \in H_*^p$ such that \begin{equation}\label{51}
\langle f, u\rangle= \sum_{i=0}^n \langle D_g^if,h_i\rangle\,\qquad
(f \in A^\infty). \end{equation} If $u$ admits a representation of
the form \eqref{51}, we will say that $u$ is \emph{of order $n$}. Note
that functionals of order $n$ are also of order $m$ for all $m > n$.
\par As in Section~\ref{28} we define a map $\Psi$ on $H_*^p$ and on
$A_*^\infty$ by \begin{equation*} (\Psi u)(w) = \langle e^{wg},u\rangle
. \end{equation*} We will often use $\hat u$ as a shorthand for $\Psi
u$. By Mergelyan's theorem, our assumption that $\Omega$ is a Jordan
domain implies that the polynomials in $g$ are dense in $A$ and thus
in $H^p$. Since \begin{equation*} (D^n\hat h)(0) = \langle g^n, h\rangle
\,\qquad(h \in H_*^p) \end{equation*} it follows that $\Psi$ is one-to-one
on $H_*^p$, that is, $\{e^{wg}\colon w \in C\mkern-12mu\mathrlap/\mkern
12mu\}$ spans a dense subspace of $H^p$. That the same is true in $A^\infty
$ follows via the same estimates that led to \eqref{51}, just like
in the usual proof of the density of trigonometric polynomials in $C^\infty
$. \par Write $X = \Psi H_*^p$ and $\mathscr X= \Psi A_*^\infty$. We
norm $X$ so that $\Psi$ is an isometry and we topologize $\mathscr
X$ so that $\Psi$ is a homeomorphism. It is clear that every $F \in
X$ is an entire function of exponential type. By \eqref{51} the same
is true of $\mathscr X$. When $p = 2$ and $g(U)$ is convex, $X$ coincides
(up to equivalence of norms) with the space $A^2(\nu)$ of Section~\ref
{28}.\par\begin{prop}\label{52} Fix $\lambda\in C\mkern-12mu\mathrlap
/\mkern12mu$ and let $S$ be a bounded operator on $H^p$ or on $A^\infty
$. Each of the following three statements about $S$ implies the other
two. \begin{enumerate}[label={\textup{(\textit{\roman*})}},listparindent=\parindent
,parsep=0pt]\item\label{53} $\mathscr R(S-R_\lambda) \subset C\mkern
-12mu\mathrlap/\mkern12mue^{\lambda g}$. \item\label{54} $(D_g-\lambda
I)S = I$. \item\label{55} For all $u$ with $\hat u(\lambda) = 0$, \begin{equation*}
(\Psi S^*u)(w) = \frac{\hat u(w)}{w-\lambda}\,\qquad(w \in C\mkern
-12mu\mathrlap/\mkern12mu). \end{equation*} \end{enumerate}\end{prop}
\begin{proof}[\textsc{proof}] Recall from Section~\ref{11} that the
operator $R_\lambda= T_g(I-\lambda T_g)^{-1}$ is bounded on $H^p$ and
on $A^\infty$, and that $(D_g-\lambda I)R_\lambda= I$ on both spaces.
Since $e^{\lambda g}$ spans the null space of $D_g-\lambda I$, \ref
{53} therefore implies \ref{54}, and \ref{54} implies that \begin{equation}\label
{56} Se^{wg} - \frac{e^{wg}}{w-\lambda} \in C\mkern-12mu\mathrlap/\mkern
12mue^{\lambda g}\,\qquad(w \in C\mkern-12mu\mathrlap/\mkern12mu),
\end{equation} which is the dual formulation of \ref{55}. If \eqref
{56} holds, then $(D_g-\lambda I)Se^{wg} = e^{wg}$, hence $(S-R_\lambda
)e^{wg} \in C\mkern-12mu\mathrlap/\mkern12mue^{\lambda g}$, for all
$w \in C\mkern-12mu\mathrlap/\mkern12mu$. Since the span of the $e^{wg}$
is dense, \ref{53} follows. \end{proof} \begin{cor}\label{57} \hfill
\begin{enumerate}[label={\textup{(\textit{\alph*})}},listparindent=\parindent
,parsep=0pt]\item\label{58} If $Y$ is $T_g$-invariant in $H^p$, then
$\Psi Y^\perp$ is nearly invariant in $X$. \item\label{59} If $\mathscr
Y$ is $D_g$-residual in $A^\infty$, then $\Psi\mathscr Y^\perp$ is
nearly invariant in $\mathscr X$. \end{enumerate}\end{cor} \begin{proof}[\textsc
{proof}] We need to show that $Y$ and $\mathscr Y$ are invariant under
some operators $S$ satisfying \ref{55}. For $Y$, take $S = R_\lambda
$. For $\mathscr Y$, let $S\colon\mathscr Y\to\mathscr Y$ be the inverse
of the restriction of $D_g-\lambda I$. Then $(S-R_\lambda)f \in C\mkern
-12mu\mathrlap/\mkern12mue^{\lambda g}$ for all $f \in\mathscr Y$,
and the Hahn-Banach theorem yields an extension of $S$ to $A^\infty
$ with $\mathscr R(S-R_\lambda) \subset C\mkern-12mu\mathrlap/\mkern
12mue^{\lambda g}$. \end{proof} \par One can show as in Section~\ref
{28} that $\Psi Y^\perp$ is in fact backward-shift invariant if we
take $g(a) = 0$, but this will be of little use to us. The role which
Theorem~\ref{38} played in Section~\ref{28} will here be delegated
to the following result, which is a more abstract version of Theorem~2.1
of \cite{7}. \begin{thm}\label{60} Suppose $Y \subset H^p$ is closed
and $N = \Psi Y^\perp$ is nearly invariant in $X$. Put \begin{equation*}
N_0 = \{F \in N\colon wF(w) \in X\}. \end{equation*} If $N_0$ is dense
in $N$, then $Y = \Theta H^p$ for some inner function $\Theta$. \end{thm}
\begin{proof}[\textsc{proof}] If $Df = f'$ on $X$ and $M_gf = gf$ on
$H^p$, a simple computation shows that $D\Psi= \Psi M_g^*$ on $H_*^p$.
In particular, $D$ is bounded on $X$. \par By our standing assumption
on $\textit d_g(\lambda)$ there exist operators $S_\lambda$ with $\Vert
S_\lambda\Vert= O(|\lambda|^{-\delta})$ satisfying \ref{55} of Proposition~\ref
{52}. For fixed $F \in N_0$, the terms of the differentiated Hadamard
product \eqref{5} are therefore $O(|\lambda_n|^{-1-\delta}) = O(n^{-1-\delta
})$ in norm, and it follows that $DF \in N$. Since $N_0$ is dense,
$N$ is $D$-invariant. In other words, $Y$ is $M_g$-invariant. Since
the polynomials in $g$ are dense in $A$, this implies shift invariance
and Beurling's theorem applies. \end{proof} \par To obtain information
about near invariance in $\mathscr X$ we will combine Theorem~\ref
{60} with an order reduction trick that should be compared to Proposition~3.1
of \cite{6} and Proposition~11 of \cite{3}. \begin{prop}\label{61}
Suppose $u \in A_*^\infty$ is of order $n$ with $\hat u(\lambda) =
0$. Define $v \in A_*^\infty$ by \begin{equation*} \hat v(w) = \frac
{\hat u(w)}{w-\lambda}\,\qquad(w \in C\mkern-12mu\mathrlap/\mkern12mu).
\end{equation*} Then $v$ is of order $0$ if $n = 0$ and of order $n-1$
if $n \ge1$. \end{prop} \begin{proof}[\textsc{proof}] By Proposition~\ref
{52}, $v = R_\lambda^*u$. The case $n = 0$ is clear. If $n \ge1$, \eqref
{51} yields $\alpha$ and $\beta$ of order $n-1$ with $u = \alpha+ D_g^*\beta
$. Because $D_gR_\lambda= I + \lambda R_\lambda$, we have $v = R_\lambda
^*\alpha+ (I+\lambda R_\lambda^*)\beta$, and the result follows by
induction. \end{proof} \begin{lem}\label{62} \hfill\begin{enumerate}[label={\textup
{(\textit{\alph*})}},listparindent=\parindent,parsep=0pt]\item\label
{63} If $Y \subset H^p$ is $T_g$-invariant and $\Theta$ is a nonconstant
inner function with $Y \subset\Theta H^p$, then $\sigma(\Theta) = \{a\}$.
\item\label{64} If $\mathscr Y\subset A^\infty$ is $D_g$-residual,
then $\mathscr Y$ has at most one common zero in $\bar{U}$. \item\label
{65} If $e^{\beta w} \in\mathscr X$, then $\beta= g(b)$ for some $b
\in\bar{U}$. \end{enumerate}\end{lem} \begin{proof}[\textsc{proof}]
\ref{63} Iteration of the formula \begin{equation*} \frac{\textit d}{\textit
dz}\int_0^zF(z,\zeta)\;\textit d\zeta= F(z,z) + \int_0^z\partial_zF(z,\zeta
)\;\textit d\zeta, \end{equation*} which holds for any $F$ holomorphic
in the bidisc, yields \begin{equation}\label{66} (T_g^{n+1}f)(z) =
\frac{1}{n!}\int_a^z[g(z)-g(\zeta)]^nf(\zeta)g'(\zeta)\;\textit d\zeta
\,\qquad(f \in H^p). \end{equation} \par Suppose $\sigma(\Theta)$ contains
some $b \ne a$. Then $b$ is a common zero of $T_gY$, and \eqref{66}
shows that $\int_a^bg^nfg' = 0$ for $f \in Y$, $n = 0,1,2,\ldots$.
Since the polynomials in $g$ are dense in $A$, the ordinary Weierstrass
approximation theorem applied on $[a,b]$ forces $Y \cap A = 0$. In
particular, $T_gY = 0$, contradicting the assumption that $Y \ne0$.
\par\ref{64} If $a$ is a common zero for $\mathscr Y$, then the $T_g$
corresponding to $a$ maps $\mathscr Y$ into itself, since $T_gD_g =
I$ on $\mathscr Y$ and $D_g\mathscr Y= \mathscr Y$. By the proof of
\ref{63}, $\mathscr Y$ can have no other common zero. \par\ref{65}
It is a consequence of Runge's theorem that if $\beta\notin g(\bar
{U})$ there exist functions $f_n$ in the span of $\{e^{wz}\colon w
\in C\mkern-12mu\mathrlap/\mkern12mu\}$ which have $f_n(\beta) = 1$
but which tend to $0$ uniformly in some neighborhood of $g(\bar{U})$.
If $\varphi_n = f_n \circ g$ and $\hat u(w) = e^{\beta w}$, it follows
that $\langle\varphi_n,u\rangle= 1$ but $\varphi_n \to0$ in $A^\infty
$, a contradiction. \end{proof} \par\begin{thm}\label{67} Suppose $1
< p < \infty$ and $g$ is a conformal map of $U$ onto a Jordan domain
with $\log g' \in A$ and $\textit d_g(\lambda) = O(|\lambda|^{-\delta
})$. Let $Y \subset H^p$ be $T_g$-invariant. \begin{enumerate}[label={\textup
{(\textit{\alph*})}},listparindent=\parindent,parsep=0pt]\item\label
{68} If $a \in U$, then $Y = B_a^nH^p$ for some $n = 1,2,3,\ldots$.
\item\label{69} If $a \in T$, then $Y = S_a^tH^p$ for some $t > 0$.
\end{enumerate}\end{thm} \textit{Remark:} We proved in Section~\ref
{28} that the backward shift on the space $A^2(\nu)$ associated to
a bounded convex domain $\Omega$ by the Lutsenko-Yulmukhametov theorem
is unicellular if $0 \in\Omega$. Theorem~\ref{67} gives the same conclusion
when $0 \in\partial\Omega$, with an uncountable chain of invariant
subspaces rather than a countable one. \begin{proof}[\textsc{proof}]
Put $M = Y^\perp$, $N = \Psi M$, and let $N_0$ be as in Theorem~\ref
{60}. By Corollary~\ref{57}, $N$ is nearly invariant. If $F \in N$
and $F(\lambda) = 0$, we thus have $F(w) = (w-\lambda)G(w)$ for some
$G \in N_0$. \par Suppose some $F \ne0$ in $N$ has a zero. Then $N_0
\ne0$ and $\bar N_0$ is nearly invariant. Theorem~\ref{60} yields $\Theta
_0$ with $\bar N_0 = \Psi Y_0^\perp$, $Y_0 = \Theta_0H^p$. Because
$Y \subset Y_0$, \ref{63} of Lemma~\ref{62} shows that $Y \cap A$ has
a common zero at $a$. Therefore, if $f \in H^p$ and $c \in C\mkern
-12mu\mathrlap/\mkern12mu$, then $T_gf + c \in Y$ only if $c = 0$.
In dual terms, $E = \{T_g^*h\colon h \in M,\,\hat h(0) = 0\}$ is dense
in $T_g^*M$. Since \begin{equation}\label{70} (\Psi T_g^*h)(w) = \frac
{\hat h(w)-e^{g(a)w}\hat h(0)}{w}\,\qquad(h \in H_*^p), \end{equation}
every $F \in\Psi E$ has $wF(w) \in X$. In other words, the annihilator
of $T_g^*M$ satisfies the hypotheses of Theorem~\ref{60}. By the Hahn-Banach
theorem, the resulting $\Theta$ has \begin{equation}\label{71} T_g(\Theta
H^p) \subset Y \subset\Theta H^p. \end{equation} \par If no $F \ne
0$ in $N$ has a zero, then $N = C\mkern-12mu\mathrlap/\mkern12mue^{\alpha
w}$ for some $\alpha\in C\mkern-12mu\mathrlap/\mkern12mu$, by Hadamard's
theorem. Since $M$ is invariant under $T_g^*$, \eqref{70} shows that
$\alpha= g(a)$ and $T_g^*M = 0$. Thus, \eqref{71} holds with $\Theta
= 1$. \par By \ref{63} of Lemma~\ref{62}, \eqref{71} forces $\Theta
= B_a^k$ if $a \in U$ and $\Theta= S_a^t$ if $a \in T$, up to unimodular
constants. It follows, as in the proofs of Propositions~\ref{46} and~\ref
{47}, that $Y$ equals $B_a^nH^p$ or $S_a^tH^p$, where $n$ is either
$k$ or $k+1$. \end{proof} \begin{thm}\label{72} Suppose $1 < p < \infty
$ and $g$ is a conformal map of $U$ onto a Jordan domain with $\log
g' \in A^\infty$ and $\textit d_g(\lambda) = O(|\lambda|^{-\delta})$
with respect to $H^p$. Let $\mathscr Y\subset A^\infty$ be $D_g$-residual.
Then $\mathscr Y= S_a^tA_a^\infty$ for some $a \in T$ and some $t \ge
0$. \end{thm} \begin{proof}[\textsc{proof}] Let $Y$ be the $H^p$-closure
of $\mathscr Y$. Put \begin{equation*} \mathscr M= \mathscr Y^\perp
\subset A_*^\infty,\quad\mathscr N= \Psi\mathscr M,\quad M = Y^\perp
\subset H_*^p,\quad N = \Psi M, \end{equation*} and let $N_0$ be as
in Theorem~\ref{60}. Note that $M$ consists precisely of those $u \in
\mathscr M$ which are of order $0$. By Corollary~\ref{57}, $\mathscr
N$ is nearly invariant in $\mathscr X$. \par Suppose some $F \ne0$
in $\mathscr N$ has infinitely many zeros. Then $N_0 \ne0$ by Proposition~\ref
{61}, and $\bar N_0$ is nearly invariant in $X$. Theorem~\ref{60} yields
$\Theta_0$ with $\bar N_0 = \Psi Y_0^\perp$, $Y_0 = \Theta_0H^p$. Because
$\mathscr Y\subset Y_0$ and $\sigma(\Theta_0) \ne\varnothing$, $\mathscr
Y$ has a common zero $a \in\bar{U}$. As in the proof of \ref{64} of
Lemma~\ref{62}, the $T_g$ corresponding to $a$ maps $Y$ into itself.
Since $\mathscr Y$ is $D_g$-invariant we must have $a \in T$, and thus
$Y = S_a^tH^p$ by Theorem~\ref{67}. By Lemma~\ref{23} it follows that
$\mathscr Y\subset S_a^tA_a^\infty$. To prove equality, we need to
show that $\langle f,u\rangle= 0$ whenever $f \in S_a^tA_a^\infty$
and $u \in\mathscr M$. \par If $\hat u$ has infinitely many zeros,
Proposition~\ref{61} shows that $\hat u = P\hat h$ for some polynomial
$P$ and some $h \in M$. Since \begin{equation*} (\Psi D_g^*u)(w) =
w\hat u(w)\,\qquad(u \in A_*^\infty) \end{equation*} and $P(D_g)f \in
S_a^tA_a^\infty\subset Y$, it follows that $\langle f, u\rangle= \langle
P(D_g)f, h\rangle= 0$. \par If $\hat u$ has only finitely many zeros,
then $\hat u(w) = P(w)e^{\beta w}$ for some polynomial $P$ and some
$\beta\in C\mkern-12mu\mathrlap/\mkern12mu$. Since $\mathscr N$ is
nearly invariant, $e^{\beta w} \in\mathscr N$. By \ref{65} of Lemma~\ref
{62}, $\beta= g(b)$ for some $b \in\bar{U}$. If $\delta_b$ is the evaluation
functional at $b$, it follows that $u = P(D_g^*)\delta_b$ and that
$b$ is a common zero for $\mathscr Y$. By \ref{64} of Lemma~\ref{62},
$b = a$, and thus $\langle f,u\rangle= 0$. \par Finally, if no $F \ne
0$ in $\mathscr N$ has infinitely many zeros, the argument of the preceding
paragraph shows that $\mathscr N= \{P(w)e^{g(a)w}\colon P\text{ polynomial}\}$
for some $a \in T$. In other words, $\mathscr Y= A_a^\infty$. \end{proof}
\section{An Integral Estimate}\label{73} Let $1 < p < \infty$ and let
$g$ be a conformal map of $U$ onto a Jordan domain $\Omega$ whose boundary
is $C^2$ with positive curvature. Recall that such $\Omega$ are always
strictly convex. \par The notation $f \lesssim h$ will mean that $f
\le Ch$ for some $C < \infty$ depending only on $g$ and $p$. If $f
\lesssim h$ and $f \gtrsim h$, we will write $f \sim h$. \par We fix
a positively oriented parametrization $\Gamma$ of $\partial\Omega$
with $|\Gamma'| = 1$; $\Gamma$ is periodic with period $L$, the length
of $\partial\Omega$. By hypothesis, $|\Gamma''| \sim1$. It will be
convenient to use the notation \begin{equation*} \Gamma_r(s) = [\Gamma
(r+s)-\Gamma(r)]\bar\Gamma'(r). \end{equation*} Note that $\Gamma_r(0)
= 0$, $\Gamma_r'(0) = 1$, and $\Gamma_r''(0) = i|\Gamma''(r)|$. \begin{lem}\label
{74} If $\Gamma_r = x_r + iy_r$ and \begin{equation*} \alpha_r(s) =
x_r(s)y_r(s/2)-y_r(s)x_r(s/2), \end{equation*} then $y_r(s) \sim|\Gamma
_r(s)|^2$ and $|\alpha_r(s)| \gtrsim|\Gamma_r(s/2)||\Gamma_r(s)|^2$.
\end{lem} \begin{proof}[\textsc{proof}] The geometry of $\Omega$ guarantees
that $|\Gamma_r|$, $y_r$, and $|\alpha_r|$ are all strictly positive
on $(0,L)$ and on $(L, 2L)$. For $|\Gamma_r|$ this is clear, for $y_r$
one may use the fact that $\partial\Omega$ intersects each of its tangents
only once, and for $|\alpha_r|$ it suffices to note that no line can
intersect $\partial\Omega$ thrice. The required estimates now follow
easily from the Taylor expansion \begin{equation*} \Gamma_r(s) = s
+ \frac{s^2}{2}[x_r''(\xi) + iy_r''(\eta)],\,\qquad|\xi|,|\eta| < |s|.\qedhere
\end{equation*} \end{proof} \begin{lem}\label{75} If $w_\lambda\in
\partial\Omega$ is the point that maximizes $\operatorname{Re}\lambda
w$, then \begin{equation*} \frac{\operatorname{Re}\lambda(w_\lambda
-z)}{|\lambda|} \sim|z-w_\lambda|^2\,\qquad(z \in\partial\Omega,\,\lambda
\ne0). \end{equation*} \end{lem} \begin{proof}[\textsc{proof}] The
left side is the distance from $z$ to the tangent of $\partial\Omega
$ at $w_\lambda$. If $w_\lambda= \Gamma(r)$ and $z = \Gamma(r+s)$,
this equals $y_r(s) \sim|\Gamma_r(s)|^2 = |z-w_\lambda|^2$, by Lemma~\ref
{74}. \end{proof} \begin{lem}\label{76} If $\zeta_t=(1-t)w+tz$, then
\begin{equation*} 1-|g^{-1}(\zeta_t)| \gtrsim t(1-t)|z-w|^2\,\qquad
(z,w \in\partial\Omega,\,0 \le t \le1). \end{equation*} \end{lem} \begin{proof}[\textsc
{proof}] By the Koebe distortion theorem (see Corollary 1.4 of \cite
{20}) the left side is comparable to the distance from $\zeta_t$ to
$\partial\Omega$. Let $a$ and $b$ be the midpoints of the complementary
arcs of $\{z,w\}$, so that $a = \Gamma(r+s/2)$ if $w = \Gamma(r)$ and
$z = \Gamma(r+s)$, $0 < s < L$. Since $\Omega$ is convex, the triangles
$[w,a,z]$ and $[z,b,w]$ lie in $\bar\Omega$. If $E$ is the union of
the four line segments $[w, a]$, $[a, z]$, $[z, b]$, $[b, w]$ we therefore
have $\textit d(\zeta_t, \partial\Omega) \ge\textit d(\zeta_t, E)$.
By symmetry it suffices to estimate $\textit d(\zeta_t,[w,a]) = \textit
d(t\Gamma_r(s),[0,\Gamma_r(s/2)])$. In the notation of Lemma \ref{74},
the line through $0$ and $\Gamma_r(s/2)$ has normal $\beta_r(s) = y_r(s/2)-ix_r(s/2)$,
hence its distance from $t\Gamma_r(s)$ is \begin{equation*} \frac{|\operatorname
{Re}t\Gamma_r(s)\bar\beta_r(s)|}{|\beta_r(s)|} = \frac{t|\alpha_r(s)|}{|\Gamma
_r(s/2)|} \gtrsim t|\Gamma_r(s)|^2 = t|z-w|^2.\qedhere\end{equation*}
\end{proof} \begin{thm}\label{77} If $1 < p < \infty$ and $g$ is a
conformal map of $U$ onto a Jordan domain whose boundary is $C^2$ with
positive curvature, then $\textit d_g(\lambda) \lesssim|\lambda|^{-\delta
}$ for some $\delta> 0$. \end{thm} \begin{proof}[\textsc{proof}] Let
$w_\lambda$ be as in Lemma~\ref{75} and define $S_\lambda$ on $H^p$
by \begin{equation*} (S_\lambda f)(\psi(z)) = \int_{w_\lambda}^{z}e^{\lambda
(z-\zeta)}f(\psi(\zeta))\;\textit d\zeta,\,\qquad(z \in\bar\Omega),
\end{equation*} where $\psi= g^{-1}$. Direct computation yields $(D_g-\lambda
I)S_\lambda= I$, so that $\textit d_g(\lambda) \le\Vert S_\lambda\Vert
$ by Proposition~\ref{52}. \par Let $z \in\partial\Omega$, let $f$
be a polynomial in $g$, and put $\zeta_t = (1-t)w_\lambda+tz$. We have
\begin{equation*} |(S_\lambda f)(\psi(z))| \le|z-w_\lambda|\bigg\{\int
_0^1e^{r\operatorname{Re}\lambda(z-\zeta_t)}\;\textit dt\bigg\}^{1/r}\bigg
\{\int_0^1|f(\psi(\zeta_t))|^{s}\;\textit dt\bigg\}^{1/s} \end{equation*}
for $1 < r < \infty$, $1/r + 1/s = 1$. By Lemma~\ref{75}, \begin{equation*}
\int_0^1e^{r\operatorname{Re}\lambda(z-\zeta_t)}\;\textit dt = \frac
{1-e^{-r\operatorname{Re}\lambda(w_\lambda-z)}}{r\operatorname{Re}\lambda
(w_\lambda-z)} \lesssim\frac{1}{r|\lambda||z-w_\lambda|^2}. \end{equation*}
By Lemma~\ref{76} and the fact that $|f(z)| \lesssim(1-|z|)^{-1/p}\Vert
f\Vert_p$ for $f \in H^p$, \begin{equation*} \bigg\{\int_0^1 |f(\psi
(\zeta_t))|^{s}\;\textit dt\bigg\}^{1/s} \lesssim C(s,p)|z-w_\lambda
|^{-2/p}\Vert f\Vert_p, \end{equation*} where \begin{equation*} C(s,p)
= \bigg\{\int_0^1t^{-s/p}(1-t)^{-s/p}\;\textit dt\bigg\}^{1/s}. \end{equation*}
In summary, for $z \in T$, \begin{equation*} |(S_\lambda f)(z)| \lesssim
C(s,p)|r\lambda|^{-1/r}|g(z)-w_\lambda|^{1-2/r-2/p}\Vert f\Vert_p.
\end{equation*} If $r$ is large enough, then $C(s,p) < \infty$ and
$|g-w_\lambda|^{1-2/r-2/p} \in L^{\mkern-2.6mup}$, uniformly in $\lambda
$. Since the polynomials in $g$ are dense in $H^p$, we are done. \end{proof}
In view of Proposition~\ref{47}, Theorems~A.2 and~B now follow from
Theorems~\ref{67} and~\ref{72}, respectively. \end{document}